\documentclass[reqno, 12pt]{amsart}
\usepackage{dsfont}
\usepackage{amssymb}
\usepackage{mathrsfs}
\usepackage{enumitem}
\usepackage{comment}
\usepackage{tikz}
\usetikzlibrary{arrows.meta}
\usepackage[backend=biber, style=alphabetic, url=false]{biblatex}
\usepackage{xurl}
\usepackage[pdfstartview=FitH, pdfborder={0 0 0}, colorlinks=true, citecolor=blue, linkcolor=blue]{hyperref}
\usepackage{aliascnt}%must load AFTER hyerref
\theoremstyle{plain}
\newtheorem{theorem}{Theorem}[section]

\newaliascnt{conjecture}{theorem}

\aliascntresetthe{conjecture}

\newaliascnt{corollary}{theorem}

\aliascntresetthe{corollary}

\newaliascnt{lemma}{theorem}
\newtheorem{lemma}[lemma]{Lemma}
\aliascntresetthe{lemma}

\newaliascnt{fact}{theorem}

\aliascntresetthe{fact}

\newaliascnt{claim}{theorem}

\aliascntresetthe{claim}

\newaliascnt{proposition}{theorem}
\newtheorem{proposition}[proposition]{Proposition}
\aliascntresetthe{proposition}

\theoremstyle{definition}

\newaliascnt{definition}{theorem}

\aliascntresetthe{definition}

\newaliascnt{example}{theorem}

\aliascntresetthe{example}

\theoremstyle{remark}

\newaliascnt{remark}{theorem}

\aliascntresetthe{remark}

\numberwithin{equation}{section}

\begin{document}
\date{}
\title[Rigidity of the parallel postulate]{The parallel postulate characterizes the Euclidean plane}

\author{Jian Chen}
\address{School of Mathematical Sciences, Laboratory of Mathematics and Complex Systems, Beijing Normal University, Beijing 100875, P. R. China}
\email{jch@mail.bnu.edu.cn}

\author{Jian Ge}
\address{School of Mathematical Sciences, Laboratory of Mathematics and Complex Systems, Beijing Normal University, Beijing 100875, P. R. China}
\email{jge@bnu.edu.cn}
\thanks{NSFC 12371049 and the Fundamental Research Funds for the Central Universities.}

\begin{abstract}
	We prove that every complete Riemannian plane satisfying the parallel axiom is isometric to the Euclidean plane. More generally, we characterize the Euclidean plane among
complete Riemannian planes without conjugate points by an asymptotic condition on geodesic intersection probabilities. This also answers affirmatively a question of Croke asking
whether a complete Riemannian plane must be Euclidean if almost every pair of geodesics intersects exactly once.
\end{abstract}
\maketitle

\section{Introduction}\label{sec:introduction}

A Riemannian plane is a smooth complete Riemannian surface diffeomorphic to $\mathbb R^2$. We say that it satisfies the \emph{parallel axiom} if, for every complete geodesic $\ell$ and every point $x\notin\ell$, there is exactly one unoriented complete geodesic through $x$ disjoint from $\ell$. This is Playfair's formulation of the parallel postulate.

The question whether this incidence condition characterizes the Euclidean metric appears in the work of Burns and Knieper \cite{BK1991}, Croke \cite[Question~1.4]{Croke2008}, and Bangert and Emmerich \cite{BE2013}. Ge, Guijarro, and Sol\'orzano \cite{GGS2020} proved the characterization for planes admitting total curvature. Our main result removes this hypothesis.

\begin{theorem}\label{thm:parallel}
	A Riemannian plane satisfying the parallel axiom is isometric to the Euclidean plane.
\end{theorem}

We prove \autoref{thm:parallel} through a stronger rigidity statement about the frequency of geodesic intersections. The connection comes from the two parts of the parallel axiom. Its existence clause implies the absence of conjugate points \cite[Proposition~1.7]{BK1991}, so every complete geodesic is a line, that is, an isometric embedding of $\mathbb R$, and distinct geodesic images meet at most once. Its uniqueness clause makes disjointness exceptional in measure: for each fixed geodesic, the directions of disjoint geodesics have angular measure zero at every point off that geodesic. Thus the parallel axiom leads naturally to an almost-everywhere intersection condition.

This condition is the subject of a question recorded by Croke \cite[Question~1.5]{Croke2008}: must a complete Riemannian plane be Euclidean if almost every pair of geodesics intersects exactly once? Croke notes that Knieper also posed this question. Our second theorem answers it affirmatively and replaces the almost-everywhere hypothesis by an asymptotic condition on intersection probabilities.

To state the result, let $(M,g)$ be a Riemannian plane without conjugate points, and fix $p\in M$. Write $d$ for the Riemannian distance, $B(x,r)$ for the open metric ball, and $\operatorname{Area}$ and $\operatorname{Length}$ for area and arclength. Let $\Gamma$ be the space of oriented complete geodesics modulo translation of their unit-speed parameters. Following Santal\'o, denote its natural geodesic density by $dG^*$; see \cite[Chapters~4 and~19]{Santalo2004}. For a Borel family $H\subset\Gamma$, write $m(H)=\int_H dG^*$, and set
\begin{equation*}
	\Gamma_r=\{\gamma\in\Gamma:\gamma\cap B(p,r)\ne\varnothing\},
	\qquad m_r=m(\Gamma_r).
	\end{equation*}
Let $I\subset\Gamma\times\Gamma$ be the set of pairs with a transverse intersection, and define
\begin{equation}\label{eq:intersection-probability}
	P(r)=\frac{1}{m_r^2}
	\int_{\Gamma_r}\int_{\Gamma_r}
	\mathbf{1}_I(\gamma,\sigma)\,dG^*(\sigma)\,dG^*(\gamma),
\end{equation}
where $\mathbf{1}_I$ is the indicator of $I$. The dependence on the fixed point $p$ is suppressed in $P$, $\Gamma_r$, and $m_r$. We verify in \autoref{sec:incidence} that $0<m_r<\infty$. Thus $P(r)$ is the probability that two geodesics sampled independently with measure $dG^*/m_r$ on $\Gamma_r$ intersect somewhere in $M$, possibly outside $B(p,r)$. Coincident pairs have measure zero, and constant rescaling of $dG^*$ does not change $P(r)$. Santal\'o considered the complementary nonintersection probability in the hyperbolic plane \cite[Section~1(c), equation~(1.6)]{Santalo1980}.

\begin{theorem}\label{thm:intersection}
	Let $M$ be a Riemannian plane without conjugate points. For any fixed $p\in M$, $\limsup_{r\to\infty}P(r)=1$ if and only if $M$ is isometric to the Euclidean plane.
\end{theorem}

Under the parallel axiom, fix $\ell$ and let $D_\ell\subset\Gamma$ be the Borel family of geodesics disjoint from $\ell$. The zero angular measure assertion above and Liouville disintegration give, for each positive integer $k$,
\begin{equation*}
	\int_{D_\ell}\operatorname{Length}(\gamma\cap B(p,k))\,dG^*(\gamma)=0.
\end{equation*}
Every line has positive length in some such ball, so $m(D_\ell)=0$. Fubini's theorem gives $P(r)=1$ for every $r>0$, and \autoref{thm:intersection} implies \autoref{thm:parallel}.

To answer Croke's question, his almost-everywhere hypothesis first excludes conjugate points: a conjugate pair would produce nearby geodesics with two transverse intersections \cite[proof of Corollary~1.3]{Croke2008}. These intersections persist under independent perturbations, and the geodesic density on local flow transversals is smooth and positive, so such pairs have positive product measure. The hypothesis then gives $P(r)=1$ at every scale, and \autoref{thm:intersection} applies.

Now we sketch the proof. The proof relates geodesic intersections to area growth through the geometry at infinity. In the ideal-boundary argument of Ge, Guijarro, and Sol\'orzano \cite[Section~3]{GGS2020}, total curvature provides the ideal boundary with its Tits metric. Uniqueness of parallels forces this boundary to be a circle of length $2\pi$, and the resulting Euclidean area growth gives flatness by \cite[Theorem~1]{BE2013}. Here square-integrability of a Riccati solution supplies the control at infinity.

Our starting point is \autoref{prop:probability-area}: there is a universal constant $C$,
\begin{equation*}
	P(r)\frac{\operatorname{Area}(B(p,r))}{r^2}\le C.
\end{equation*}
Thus $\limsup_{r\to\infty}P(r)>0$ already gives bounded normalized area along a sequence. Let $U$ be the Riccati solution supplied by \cite[Lemma~2.2]{BE2013}, and let $E$ be the integral of $U^2$ with respect to Liouville measure, divided by $2\pi$. The estimate in \cite[proof of Theorem~1, p.~374]{BE2013} gives $E<\infty$; see \autoref{prop:integrability}. We show that this integrability gives local Euclidean estimates on balls escaping to infinity and convergence of the rescaled planes to a unique Euclidean cone. By \autoref{prop:area-asymptotics} and \autoref{thm:asymptotic-cone}, its total angle $\Theta$ and the area growth satisfy
\begin{equation*}
\Theta=2\pi+E,\qquad
\lim_{r\to\infty}\frac{\operatorname{Area}(B(p,r))}{r^2}=\pi+\frac E2.
\end{equation*}

If $E>0$, the excess angle allows two separated families of complete cone lines. A four-point consequence of Crofton's formula gives the corresponding families on $M$ positive normalized measure; see \autoref{sec:separated-families}. Cone convergence excludes intersections between these families in any fixed rescaled ball. The tail estimate in \autoref{prop:intersection-tail} controls intersections that might escape to infinity. Together these imply $\limsup_{r\to\infty}P(r)<1$. Thus the hypothesis of \autoref{thm:intersection} forces $E=0$ and Euclidean area growth, and the equality case of \cite[Theorem~1]{BE2013} completes the proof.

\section{Incidence and a quadratic area bound}\label{sec:incidence}
Throughout the rest of the paper, $(M,g)$ is a Riemannian plane without conjugate points, and $p$ is the fixed point chosen above. Let $T^1M$ be its unit tangent bundle and $T^1_xM$ the unit tangent circle at $x\in M$. We use Santal\'o's notation $dP$ for the area element and $d\theta_x$ for angular measure of mass $2\pi$ on $T^1_xM$.
The exponential map at every point is a diffeomorphism. Consequently every complete geodesic minimizes between any two of its points and is a proper embedding. Two distinct geodesic images meet at most once, and their intersection, if present, is transverse; compare \cite[Lemma 1.4]{GG2023}. We denote the unique closed geodesic segment joining $a,b\in M$ by $[a,b]$.

For $v\in T^1M$, let $\gamma_v$ be the complete unit-speed geodesic with $\dot\gamma_v(0)=v$. Let $\pi:T^1M\to M$ be the projection and $\Phi_t$ the geodesic flow. The flow action is free and proper, since
\begin{equation*}
d(\pi v,\pi\Phi_t v)=|t|.
\end{equation*}
Thus $\Gamma=T^1M/\mathbb R$ is a smooth surface. The geodesic density is the oriented lift of Santal\'o's density $dG$; forgetting orientation pushes $dG^*$ forward to $2dG$. In local flow coordinates its normalization is
\begin{equation}\label{eq:santalo-disintegration}
dG^*\,dt=dP\,d\theta_x.
\end{equation}
This is the oriented convention in \cite[equation~(4.55)]{Santalo2004}, together with the Riemannian disintegration in \cite[Sections~19.2, 19.4--19.5]{Santalo2004}. Products here denote positive densities. In coordinates given by arclength $s$ on a transversal and crossing angle $\beta\in[0,2\pi)$,
\begin{equation*}
dG^*=|\sin\beta|\,ds\,d\beta.
\end{equation*}
Consequently, Crofton's formula gives
\begin{equation}\label{eq:crofton-segment}
m\{\gamma:\gamma\text{ crosses }[a,b]\}=4d(a,b).
\end{equation}
The pencils through the endpoints and the lines containing a given segment have measure zero.

For a Borel family $H\subset\Gamma$, put
\begin{equation*}
I_x(H)=\{v\in T^1_xM:\gamma_v\in H\}.
\end{equation*}
If $D\subset M$ is Borel, let $N_H(D)$ be the product measure of ordered pairs in $H^2$ whose transverse intersection lies in $D$. Each pair has only one intersection. In the first integral below, $x=\gamma(t)$. The pair density is obtained from the transversal formula and \eqref{eq:santalo-disintegration}, as in \cite[equation~(4.19)]{Santalo2004}:
\begin{equation}\label{eq:intersection-disintegration}
\begin{aligned}
N_H(D)
&=\int_H\int_{\gamma\cap D}\int_{I_x(H)}
|\sin\angle(\dot\gamma,w)|\,d\theta_x(w)\,dt\,dG^*(\gamma)\\
&=\int_D\int_{I_x(H)}\int_{I_x(H)}
|\sin\angle(v,w)|\,d\theta_x(w)\,d\theta_x(v)\,dP.
\end{aligned}
\end{equation}
All these identities concern oriented lines and ordered pairs; no additional factor is introduced. Tonelli's theorem applies before finiteness has been established.

A line spends at most $2r$ units of time in $B(p, r)$, and every line meeting $B(p, r)$ spends at least $2r$ units of time in $B(p,2r)$. Therefore
\begin{equation}\label{eq:line-mass-area}
	2\pi \operatorname{Area}(B(p,r))\le 2rm_r\le 2\pi \operatorname{Area}(B(p,2r)).
\end{equation}
In particular $0<m_r<\infty$.

The main estimate of this section is the following. Only its finiteness, rather than the numerical constant, will be used.

\begin{proposition}\label{prop:probability-area}
	For every $r>0$,
	\begin{equation}\label{eq:probability-area}
		P(r)\frac{\operatorname{Area}(B(p,r))}{r^2}\le32.
	\end{equation}
	Consequently, if $\limsup_{r\to\infty}P(r)>0$, there is a sequence $r_i\to\infty$ with bounded ratios $\operatorname{Area}(B(p,r_i))/r_i^2$.
\end{proposition}

Two elementary estimates enter the proof. The first converts the angular intersection density into a quadratic radial error.
\begin{lemma}\label{lem:angular-integral}
	For any measurable $I\subset\mathbb R/(2\pi\mathbb Z)$ and any $\alpha$,
	\begin{equation*}
		\int_I\int_I|\sin(\theta-\varphi)|\,d\theta\,d\varphi	\le 8\pi\int_I\sin^2(\theta-\alpha)\,d\theta.
	\end{equation*}
\end{lemma}
\begin{proof}
	Put $f(\theta)=|\sin(\theta-\alpha)|$, and write $|J|$ for the angular measure of a measurable set $J$. For $0\le t\le1$, we have $|\{f\le t\}|=4\arcsin t\le2\pi t$. On $f(\varphi)\le f(\theta)$, the sine difference formula gives $|\sin(\theta-\varphi)|\le2f(\theta)$. By symmetry,
	\begin{align*}
		\int_{I^2}|\sin(\theta-\varphi)|\,d\theta\,d\varphi
		&\le4\int_I f(\theta)\,|I\cap\{f\le f(\theta)\}|\,d\theta\\
		&\le8\pi\int_I f(\theta)^2\,d\theta.
	\end{align*}
\end{proof}

\begin{lemma}\label{lem:two-ended-excess}
	Let $\gamma:\mathbb R\to M$ be a unit-speed complete geodesic and let $h(t)=d(p,\gamma(t))$. For every $s\ge0$,
	\begin{equation}\label{eq:two-ended-excess}
		\int_{|t|\ge s}(1-h'(t)^2)\,dt
		\le2\bigl(h(s)+h(-s)-2s\bigr).
	\end{equation}
\end{lemma}
\begin{proof}
	The functions $\delta_\pm(t)=h(\pm t)-t$, $t\ge0$, decrease and are bounded below by $-h(0)$. Their limits $a_\pm$ satisfy $a_++a_-\ge0$, since $h(t)+h(-t)\ge2t$. Applying $1-q^2\le2(1-q)$ to the outgoing radial derivative on each half-line yields
	\begin{equation*}
		\int_{|t|\ge s}(1-h'(t)^2)\,dt
		\le2\bigl(\delta_+(s)-a_++\delta_-(s)-a_-\bigr),
	\end{equation*}
	which proves the assertion.
\end{proof}

\begin{proof}[Proof of \autoref{prop:probability-area}]
	At $x\ne p$, choose the radial direction $\nabla d(p,\cdot)$ as $\alpha$ in \autoref{lem:angular-integral}. Along $\gamma$, put $h(t)=d(p,\gamma(t))$. The scalar product of the radial direction with $\dot\gamma$ is $h'$. Equations \eqref{eq:intersection-disintegration} and \eqref{eq:two-ended-excess}, with the origin of the parameter chosen in $B(p,r)$, imply
	\begin{equation*}
		P(r)m_r^2
		\le C_0\int_{\Gamma_r}\int_{\mathbb R}(1-h'^2)\,dt\,dG^*
		\le4C_0rm_r.
	\end{equation*}
	The inner integral is unchanged by translating the parameter, so the estimate holds on line space independently of this choice. The lower bound in \eqref{eq:line-mass-area} now gives \eqref{eq:probability-area}.
\end{proof}

\section{The Riccati equation and area growth}\label{sec:riccati}

Let $K$ be the Gaussian curvature and $X$ the geodesic vector field. Hopf's construction \cite{Hopf1948} gives a Borel measurable, locally bounded function $U:T^1M\to\mathbb R$, differentiable along flow lines, satisfying
\begin{equation*}
XU+U^2+K\circ\pi=0;
\end{equation*}
see \cite[Lemma~2.2]{BE2013}. With $d\omega_x=d\theta_x/(2\pi)$, set
\begin{equation*}
e(x)=\int_{T^1_xM}U(v)^2\,d\omega_x(v),\qquad
Z(x)=\int_{T^1_xM}U(v)v\,d\omega_x(v),
\end{equation*}
and $E=\int_M e\,dP$.

For $\zeta\in C_c^\infty(M)$, flow invariance of Liouville measure applied to $(\zeta\circ\pi)U$ gives
\begin{equation}\label{eq:averaged-riccati}
\int_M(K+e)\zeta\,dP=\int_M\langle Z,\nabla\zeta\rangle\,dP,
\qquad |Z|^2\le e.
\end{equation}
By local boundedness, the identity extends to compactly supported Lipschitz test functions.

Constants $C$ are universal, whereas $C_E$ depend only on $E$; their values may change. Write $a_+=\max\{a,0\}$.

\begin{proposition}\label{prop:integrability}
	If $\operatorname{Area}(B(p,r_i))/r_i^2$ is bounded along some sequence $r_i\to\infty$, then $E<\infty$. Conversely, finite $E$ implies
	\begin{equation}\label{eq:uniform-area-bounds}
		2r^2\le \operatorname{Area}(B(x,r))\le C_Er^2\qquad(x\in M,\ r>0).
	\end{equation}
\end{proposition}

\begin{proof}
	The first assertion follows from \cite[proof of Theorem~1, p.~374]{BE2013}, whose estimate in our normalization is
	\begin{equation*}
		E\le2\left(\liminf_{r\to\infty}
		\frac{\operatorname{Area}(B(p,r))}{r^2}-\pi\right).
	\end{equation*}
	For the converse, suppose $E<\infty$ and put $\psi_{x,r}=(1-d(x,\cdot)/r)_+$. Coarea, first variation, Gauss--Bonnet, and \eqref{eq:averaged-riccati} give
	\begin{equation}\label{eq:linear-radial-cutoff}
		\frac{\operatorname{Length}(\partial B(x,r))}{r}
		=2\pi+\int_M e\psi_{x,r}\,dP
		-\int_M\langle Z,\nabla\psi_{x,r}\rangle\,dP.
	\end{equation}
	Multiply by the radius and integrate. By Cauchy--Schwarz and monotonicity of area, $a=\operatorname{Area}(B(x,r))/r^2$ satisfies
	\begin{equation*}
		\begin{aligned}
		a&\le\pi+\frac E2+\frac{\sqrt E}{r^2}
		\int_0^r\sqrt{\operatorname{Area}(B(x,t))}\,dt\\
		&\le\pi+\frac E2+\sqrt{Ea}\le\pi+E+\frac a2.
		\end{aligned}
	\end{equation*}
	Thus $a\le2\pi+2E$. A geodesic through $x$ cuts $\partial B(x,r)$ into two arcs whose endpoints are distance $2r$ apart. Hence $\operatorname{Length}(\partial B(x,r))\ge4r$; integration gives the lower bound.
\end{proof}

\begin{proposition}\label{prop:area-asymptotics}
	If $E<\infty$, then
	\begin{equation}\label{eq:fixed-center-asymptotics}
		\begin{aligned}
		\lim_{r\to\infty}\frac{\operatorname{Length}(\partial B(p,r))}r&=2\pi+E,\\
		\lim_{r\to\infty}\frac{\operatorname{Area}(B(p,r))}{r^2}&=\pi+\frac E2.
		\end{aligned}
	\end{equation}
\end{proposition}

\begin{proof}
	Put $\psi_r=\psi_{p,r}$. Dominated convergence gives $\int_M e\psi_r\,dP\to E$. For each $D>0$, \eqref{eq:uniform-area-bounds} implies
	\begin{equation*}
		\left|\int_M\langle Z,\nabla\psi_r\rangle\,dP\right|
		\le\frac1r\int_{B(p,D)}|Z|\,dP
		+\sqrt{C_E}\,\|Z\|_{L^2(M\setminus B(p,D))}.
	\end{equation*}
	Since $|Z|^2\le e$, letting $r\to\infty$ and then $D\to\infty$ makes this term vanish. Equation \eqref{eq:linear-radial-cutoff} gives the perimeter limit, and integration gives the area limit.
\end{proof}

Assume $E<\infty$. For $y\in M$ and $a>0$ with $\eta=\int_{B(y,a)}e\,dP\le1$, equations \eqref{eq:linear-radial-cutoff} and \eqref{eq:uniform-area-bounds} give, for $0<t\le a$,
\begin{equation}\label{eq:local-perimeter-estimate}
\left|\frac{\operatorname{Length}(\partial B(y,t))}t-2\pi\right|
\le\eta+\frac{\sqrt{\eta \operatorname{Area}(B(y,t))}}t\le C_E\sqrt\eta.
\end{equation}

\section{The asymptotic cone}\label{sec:asymptotic-cone}

For planes of finite total curvature, convergence to the cone over the ideal boundary is classical; see \cite[Section~3.7, Theorem~3.7.2]{SST2003}. We establish the cone structure here from the assumption $E<\infty$.

For $\Lambda>0$, let $C_\Lambda$ be the metric completion of
\begin{equation*}
(0,\infty)\times(\mathbb R/\Lambda\mathbb Z),\qquad
d\rho^2+\rho^2d\vartheta^2,
\end{equation*}
and denote its vertex by $0$. We shall prove the following.

\begin{theorem}\label{thm:asymptotic-cone}
	If $E<\infty$, put $\Theta=2\pi+E$. Then
	\begin{equation*}
		(M,r^{-1}d,p)\longrightarrow(C_\Theta,0)
	\end{equation*}
	in the pointed Gromov--Hausdorff topology as $r\to\infty$.
\end{theorem}

The uniform area bounds give metric precompactness. To identify the limits, we first turn small integrals of $e$ into Euclidean distance and area estimates. The center is allowed to move slightly so that the Riccati integral is small along most radial directions.

\begin{lemma}\label{lem:local-distance-chart}
	Assume $E<\infty$, and let $x\in M$ and $r>0$ satisfy $\eta=\int_{B(x,10r)}e\,dP\le1$. There are $y\in B(x,r/8)$ and a radius-preserving bijection
	\begin{equation*}
		F:B(y,r)\longrightarrow B_{\mathbb R^2}(0,r)
	\end{equation*}
	such that
	\begin{equation}\label{eq:local-distance-distortion}
		\sup_{z,w\in B(y,r)}\bigl|d(z,w)-|F(z)-F(w)|\bigr|
		\le C_Er\eta^{1/4}.
	\end{equation}
	The distance on the left is the ambient distance in $M$.
	Writing $dP_0$ for Euclidean area, the same map satisfies
	\begin{equation}\label{eq:local-area-distortion}
		\bigl|F_*(dP|_{B(y,r)})-dP_0\bigr|\bigl(B_{\mathbb R^2}(0,r)\bigr)
		\le C_Er^2\sqrt\eta,
	\end{equation}
	where the absolute value denotes the total variation of a signed measure.
\end{lemma}

\begin{proof}
	Under the rescaling $g\mapsto r^{-2}g$, the Riccati function scales by $r$ and the measure $e\,dP$ is unchanged. We may therefore suppose $r=1$. For $y\in B(x,1/8)$ and $v\in T^1_yM$, put
	\begin{equation*}
		D_y(v)=\int_0^8|U(\Phi_t v)|\,dt.
	\end{equation*}
	Cauchy--Schwarz and flow invariance give
	\begin{equation*}
		\int_{B(x,1/8)}\int_{T^1_yM}D_y(v)^2\,d\omega_y(v)\,dP
		\le64\eta,
	\end{equation*}
	because the flow segments involved stay in $B(x,10)$. The lower area bound in \eqref{eq:uniform-area-bounds} supplies a point $y\in B(x,1/8)$ with
	\begin{equation*}
		\int_{T^1_yM}D_y(v)^2\,d\theta_y(v)\le C\eta.
	\end{equation*}

	Parametrize $T^1_yM$ by angle, writing $v_\theta$ for the vector of angle $\theta$. Put $D_y(\theta)=D_y(v_\theta)$. In the corresponding polar coordinates about $y$, write $g=dt^2+J(t,\theta)^2d\theta^2$ and set $F(t,\theta)=t(\cos\theta,\sin\theta)$. Along the ray with initial vector $v_\theta$, the positive function
	\begin{equation*}
		j(t)=\exp\left(\int_0^t U(\Phi_u v_\theta)\,du\right)
	\end{equation*}
	solves the Jacobi equation and has $j(0)=1$. Set $a(\theta)=\exp(-3D_y(\theta))$. Reduction of order gives
	\begin{equation*}
		J(t,\theta)=j(t)\int_0^t j(u)^{-2}\,du
		\ge t a(\theta),\qquad 0\le t\le8.
	\end{equation*}
	Thus $\int_0^{2\pi}(1-a)\,d\theta\le C\sqrt\eta$. Using $|J-t|=J-t+2(t-J)_+$ and \eqref{eq:local-perimeter-estimate}, we obtain
	\begin{equation}\label{eq:jacobi-area-error}
		\int_0^{2\pi}|J(t,\theta)-t|\,d\theta
		\le C_Et\sqrt\eta,\qquad 0<t\le8.
	\end{equation}
	Since $F_*dP=(J/t)dP_0$ away from the origin, integrating this estimate over $0<t<1$ proves \eqref{eq:local-area-distortion} at the normalized scale.

	For the lower distance estimate, change the angular coordinate by
	\begin{equation*}
		H(\theta)=\frac{2\pi\int_0^\theta a(\phi)\,d\phi}
		{\int_0^{2\pi}a(\phi)\,d\phi}.
	\end{equation*}
	When $\eta$ is sufficiently small, $H$ is Lipschitz and satisfies
	\begin{equation*}
		\|H-\operatorname{id}\|_\infty\le C\sqrt\eta,
		\qquad tH'(\theta)\le(1+C\sqrt\eta)J(t,\theta)
	\end{equation*}
	almost everywhere. A shortest segment joining points of $B(y,1)$ stays in $B(y,3)$. If it is nonradial, its angular derivative never vanishes: a radial tangency would force the entire geodesic to be radial. The one-dimensional chain rule therefore applies to $H$ along this segment. Replacing the polar angle $\theta$ by $H(\theta)$ increases the length by at most the factor $1+C\sqrt\eta$. Radial pieces preserve their length. Correcting the endpoint angles consequently gives
	\begin{equation*}
		|F(z)-F(w)|\le d(z,w)+C\sqrt\eta.
	\end{equation*}
	This argument uses a minimizing segment in the whole surface, as required for the ambient distance estimate.

	For the upper estimate, a nonradial Euclidean chord in the unit disk has polar equation $t=f(\theta)$ on an interval $I$ of length at most $\pi$. Replace it by $t=f(\theta)+h$, with $0\le h\le b\le1$, and connect its endpoints radially to the original endpoints. Since $u\mapsto\sqrt{f'(\theta)^2+u^2}$ is $1$-Lipschitz for $u\ge0$, the lifted path has length at most
	\begin{equation*}
		|F(z)-F(w)|+(2+\pi)b
		+\int_I\bigl(J(f(\theta)+h,\theta)-f(\theta)-h\bigr)_+\,d\theta.
	\end{equation*}
	Averaging the last term over $0\le h\le b$ and using \eqref{eq:jacobi-area-error} bounds it, for some $h$, by
	\begin{equation*}
		\frac1b\int_0^2\int_0^{2\pi}(J(t,\theta)-t)_+\,d\theta\,dt
		\le C_E\frac{\sqrt\eta}{b}.
	\end{equation*}
	Choose $b=\eta^{1/4}$ if $\eta>0$, and let $b\downarrow0$ if $\eta=0$. Radial chords lift to radial minimizing segments and are handled directly. For $\eta$ bounded away from zero, the diameter bounds suffice after increasing $C_E$. Restoring the scale proves both estimates.
\end{proof}

\begin{proof}[Proof of \autoref{thm:asymptotic-cone}]
	Fix a sequence $r_i\to\infty$, and put $d_i=r_i^{-1}d$. We apply the metric form of Gromov's compactness theorem \cite[Theorem~8.1.10]{BBI2001}. To verify its covering hypothesis, fix $b,\delta>0$ and let $x_1,\ldots,x_N$ be $\delta$-separated points in $B_{d_i}(p,b)$. The original metric balls $B(x_j,\delta r_i/2)$ are pairwise disjoint and contained in $B(p,(b+\delta/2)r_i)$. Thus \eqref{eq:uniform-area-bounds} gives
	\begin{equation*}
		2N\left(\frac{\delta r_i}{2}\right)^2
		\le C_E\left(b+\frac\delta2\right)^2r_i^2,
		\qquad N\le2C_E\left(\frac b\delta+\frac12\right)^2.
	\end{equation*}
	A maximal $\delta$-separated set is a $\delta$-net, so the required bound on the number of net points is uniform in $i$. The compactness theorem therefore gives, after passage to a subsequence, a proper geodesic limit $(Y,d_Y,p_*)$. This application uses the area bounds and does not require a curvature lower bound. In a common realization of the convergence, the uniform area bounds also give a subsequential locally weak limit $\nu$ of $r_i^{-2}dP$.

	We first show that $Y\setminus\{p_*\}$ is locally Euclidean. Let $x_i\to x\ne p_*$ and choose $0<s<d_Y(p_*,x)/20$. The original balls $B(x_i,10sr_i)$ leave every compact subset of $M$, so the integrals of $e$ over these balls tend to zero. By \autoref{lem:local-distance-chart}, there are $y_i\in B_{d_i}(x_i,s/8)$ and radius-preserving bijections $F_i$ from $B_{d_i}(y_i,s)$ to the Euclidean $s$-ball whose distortions tend to zero. Take a subsequence with $y_i\to y$. Moving approximating points radially to any smaller radius shows convergence of the corresponding closed balls. A diagonal subsequence of the inverse maps converges uniformly on compact Euclidean subballs to an isometry onto $B_Y(y,s)$; denote its inverse by $F$. Since $d_Y(x,y)\le s/8$, this is a Euclidean neighborhood of $x$. Weak convergence and \eqref{eq:local-area-distortion}, tested against continuous functions compactly supported in this chart, show that $F$ carries $\nu|_{B_Y(y,s)}$ to Euclidean area on $B_{\mathbb R^2}(0,s)$. Hence $\nu$ is ordinary flat area on $Y\setminus\{p_*\}$. The uniform upper area bound also gives $\nu(\{p_*\})=0$.

	Next, every point $x\ne p_*$ lies on a minimizing ray from $p_*$. To obtain one, extend the radial segments through approximating points in $M$ and take a limit on compact parameter intervals. This ray is unique. In fact, concatenate any other minimizing segment from $p_*$ to $x$ with the outgoing part of the chosen ray. The concatenation still minimizes to every later point. In a Euclidean neighborhood of $x$ it cannot have a corner, so the incoming segments agree near $x$. Unique continuation in the flat punctured surface makes them agree all the way to $p_*$. The outgoing continuation is unique for the same reason. Compactness and this uniqueness imply continuous dependence of radial rays on any marked point away from $p_*$.

	Write $\rho(z)=d_Y(p_*,z)$. At $x\ne p_*$, choose nearby points $x_-$ and $x_+$ on its radial ray, before and after $x$, all in one Euclidean neighborhood. For $z$ near $x$,
	\begin{equation*}
		\rho(x_+)-d_Y(x_+,z)\le\rho(z)
		\le\rho(x_-)+d_Y(x_-,z).
	\end{equation*}
	The bounding functions are smooth near $x$ and have the same value and first derivative there. Thus $\rho$ is differentiable, with gradient the radial unit vector. Continuous dependence of the rays gives $\rho\in C^1(Y\setminus\{p_*\})$ and $|\nabla\rho|=1$. Each level $S_t=\rho^{-1}(t)$ is consequently a compact one-manifold. It is connected: the distance circles in the approximating planes converge to $S_t$ in Hausdorff distance, since any point with radius $t+o(1)$ can be moved radially to radius exactly $t$ at a cost $o(1)$. Hence $S_t$ is a circle. Radial transport now gives a homeomorphism
	\begin{equation*}
		Y\setminus\{p_*\}\simeq(0,\infty)\times S_1.
	\end{equation*}

	The circles in the approximating planes, parametrized with constant speed on $[0,1]$, have a subsequential limiting parametrization onto $S_t$. By \eqref{eq:fixed-center-asymptotics} and lower semicontinuity of length,
	\begin{equation*}
		\operatorname{Length}(S_t)\le\Theta t,\qquad \Theta=2\pi+E.
	\end{equation*}
	Put $\Lambda=\operatorname{Length}(S_1)$, and lift an arclength parameter $s\bmod\Lambda$ on $S_1$ to the universal cover of the flat annulus $Y\setminus\{p_*\}$. Under a developing map to $\mathbb R^2$, the lifted radial rays have the form
	\begin{equation*}
		\operatorname{dev}(\rho,s)=a(s)+\rho v(s),\qquad |v(s)|=1.
	\end{equation*}
	For $s,s'$ in one period, their radius-$t$ points are joined by a lifted circle arc of length at most $\Theta t$, so $|a(s)-a(s')|\le(\Theta+2)t$. Letting $t\downarrow0$, and comparing different periods by a finite chain, shows that $a(s)$ is a common point $a$.

	The developed unit level $a+v(s)$ is parametrized by arclength, so $|v'(s)|=1$ and $\langle v(s),v'(s)\rangle=0$. Thus the metric is $d\rho^2+\rho^2ds^2$, with period $\Lambda$ in $s$, identifying the punctured surface with the punctured cone. A minimizing segment through $p_*$ can be replaced inside radius $t$ by an arc of $S_t$, at an additional cost at most $\Theta t$. Letting $t\downarrow0$ identifies its metric completion with $Y=C_\Lambda$.

	Since $\nu$ is cone area, $S_1$ has zero $\nu$-measure. Weak convergence and \eqref{eq:fixed-center-asymptotics} now give
	\begin{equation*}
		\frac\Lambda2=\nu(B_Y(p_*,1))
		=\lim_i\frac{\operatorname{Area}(B(p,r_i))}{r_i^2}
		=\pi+\frac E2.
	\end{equation*}
	All subsequential limits therefore have angle $\Theta$. Metric precompactness then proves convergence of the full family.
\end{proof}

\section{Separated families of geodesics}\label{sec:separated-families}

The cone theorem reduces the remaining problem to excluding a positive angle excess. We shall prove the following quantitative form of this assertion.

\begin{proposition}\label{prop:probability-gap}
	If $0<E<\infty$, there is $\varepsilon>0$, depending only on $E$, such that
	\begin{equation*}
		\limsup_{r\to\infty}P(r)\le1-\varepsilon.
	\end{equation*}
\end{proposition}

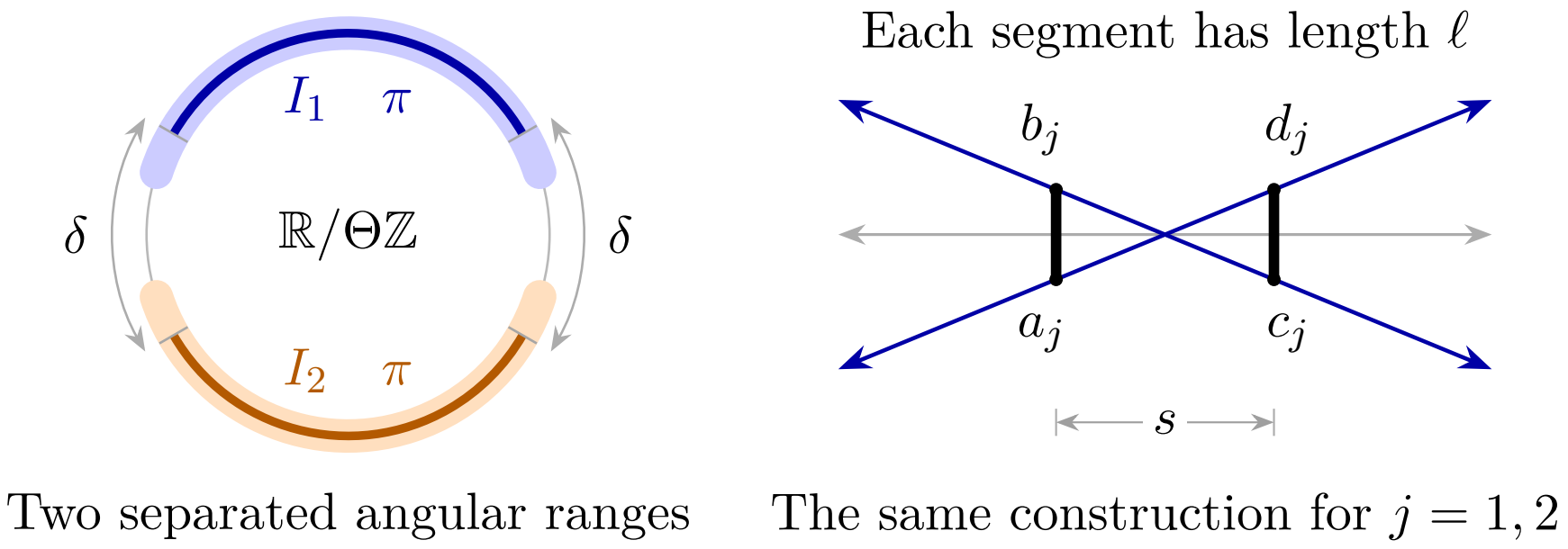
\begin{figure}[htbp]
	\centering
	\begin{tikzpicture}[x=1cm,y=1cm,font=\small,>=Stealth]
		% The circular drawing is schematic: labels indicate intrinsic arc lengths.
		\draw[gray!55,line width=.5pt] (2.65,2.55) circle[radius=1.48];
		\draw[blue!20,line width=7pt,line cap=round] (2.65,2.55) ++(18:1.48) arc[start angle=18,end angle=162,radius=1.48];
		\draw[orange!25,line width=7pt,line cap=round] (2.65,2.55) ++(198:1.48) arc[start angle=198,end angle=342,radius=1.48];
		\draw[blue!65!black,line width=1.8pt] (2.65,2.55) ++(30:1.48) arc[start angle=30,end angle=150,radius=1.48];
		\draw[orange!70!black,line width=1.8pt] (2.65,2.55) ++(210:1.48) arc[start angle=210,end angle=330,radius=1.48];
		\foreach \angle in {30,150,210,330} {
		\draw[gray!70,line width=.5pt] (2.65,2.55) ++(\angle:1.36) -- ++(\angle:.24);
		}
		\node[blue!65!black] at (2.65,3.55) {$I_1\quad\pi$};
		\node[orange!70!black] at (2.65,1.55) {$I_2\quad\pi$};
		\node at (2.65,2.55) {$\mathbb R/\Theta\mathbb Z$};
		\draw[<->,gray!65,line width=.5pt] (2.65,2.55) ++(150:1.72) arc[start angle=150,end angle=210,radius=1.72];
		\draw[<->,gray!65,line width=.5pt] (2.65,2.55) ++(330:1.72) arc[start angle=330,end angle=390,radius=1.72];
		\node at (.65,2.55) {$\delta$};
		\node at (4.65,2.55) {$\delta$};
		\node at (2.65,.48) {Two separated angular ranges};
		% A local developed view; arrows denote continuation beyond this window.
		\draw[<->,gray!65,line width=.6pt] (6.25,2.55) -- (11.05,2.55);
		\draw[<->,blue!65!black,line width=.85pt] (6.25,1.56) -- (11.05,3.54);
		\draw[<->,blue!65!black,line width=.85pt] (6.25,3.54) -- (11.05,1.56);
		\draw[line width=2.2pt] (7.85,2.22) -- (7.85,2.88);
		\draw[line width=2.2pt] (9.45,2.22) -- (9.45,2.88);
		\fill (7.85,2.22) circle[radius=1.4pt];
		\fill (7.85,2.88) circle[radius=1.4pt];
		\fill (9.45,2.22) circle[radius=1.4pt];
		\fill (9.45,2.88) circle[radius=1.4pt];
		\node[below,xshift=-3pt,yshift=-5pt,inner sep=2pt] at (7.85,2.22) {$a_j$};
		\node[above,xshift=-3pt,yshift=5pt,inner sep=2pt] at (7.85,2.88) {$b_j$};
		\node[below,xshift=3pt,yshift=-5pt,inner sep=2pt] at (9.45,2.22) {$c_j$};
		\node[above,xshift=3pt,yshift=5pt,inner sep=2pt] at (9.45,2.88) {$d_j$};
		\draw[gray!65] (7.85,1.27) -- (7.85,1.07);
		\draw[gray!65] (9.45,1.27) -- (9.45,1.07);
		\draw[<->,gray!75] (7.85,1.17) -- node[fill=white,inner sep=2pt,text=black] {$s$} (9.45,1.17);
		\node at (8.65,4.02) {Each segment has length $\ell$};
		\node at (8.65,.48) {The same construction for $j=1,2$};
	\end{tikzpicture}
	\caption{The construction in \autoref{lem:cone-families}: disjoint enlarged angular ranges (left) and transversal segments selecting each family of cone lines (right).}
	\label{fig:cone-families}
\end{figure}

We first describe the geometric configuration responsible for the gap. A complete cone line avoiding the vertex sweeps an angular interval of length $\pi$. When $\Theta>2\pi$, two such intervals fit on the circle of directions with a positive gap on either side. These gaps allow both lines to vary while keeping the two families disjoint; see \autoref{fig:cone-families}. Two short transversal segments will select each family, whose size will be transferred to $M$ using Crofton's formula. The possibility that intersections disappear to infinity will then be addressed in \autoref{sec:distant-intersections}.

\begin{lemma}\label{lem:cone-families}
	If $\Theta>2\pi$, there are two pairs of geodesic segments $S_-^j=[a_j,b_j]$ and $S_+^j=[c_j,d_j]$, $j=1,2$, lying strictly inside $B_{C_\Theta}(0,1)$, with the following properties. Every complete line meeting both segments of the first pair is disjoint from every complete line meeting both segments of the second pair. Moreover, for some $q>0$ depending only on $\Theta$,
	\begin{equation}\label{eq:cone-four-point-gap}
		2\bigl|d(a_j,c_j)+d(b_j,d_j)-d(a_j,d_j)-d(b_j,c_j)\bigr|=q,
		\qquad j=1,2.
	\end{equation}
\end{lemma}

\begin{proof}
	Put $\delta=(\Theta-2\pi)/2$. In the circle $\mathbb R/\Theta\mathbb Z$, the intervals
	\begin{equation*}
		I_1=[0,\pi],\qquad I_2=[\pi+\delta,2\pi+\delta]
	\end{equation*}
	are separated by two gaps of length $\delta$. Each corresponding closed sector is isometric, for the ambient cone distance, to a Euclidean closed half-plane. Indeed, the cone distance is given by
	\begin{equation*}
		d\bigl((r,\varphi),(s,\psi)\bigr)^2
		=r^2+s^2-2rs\cos\!\left(\min\{\pi,d_{\mathbb R/\Theta\mathbb Z}(\varphi,\psi)\}\right).
	\end{equation*}
	A minimizing segment avoiding the vertex develops to a Euclidean chord, whose continuous angular variation is less than $\pi$; otherwise the segment passes through the vertex. Since $\Theta>2\pi$, the displayed formula gives the Euclidean chord distance within each indicated sector.

	Use upper-half-plane coordinates $(x,y)$ on each sector. Put $s=1/4$, $b_0=1/2$, and choose $0<\ell<1/8$ so that $2\arctan(\ell/s)<\delta$. Write $\alpha=\arctan(\ell/s)$ and set
	\begin{equation*}
		S_\pm^j=\{(\pm s/2,y):|y-b_0|\le\ell/2\},\qquad j=1,2.
	\end{equation*}
	These segments lie strictly inside the unit ball. A complete line meeting both segments of pair $j$ contains the unique Euclidean chord joining the two marked points. Its developed affine line has equation $y=kx+b$, with
	\begin{equation*}
		|k|\le\ell/s,\qquad b\ge b_0-\ell/2>0.
	\end{equation*}
	It avoids the vertex and has a unique continuation. Its entire angular sweep lies in $I_j$ enlarged at each end by $\alpha$. The two enlarged intervals are disjoint since their intervening gaps have length $\delta-2\alpha>0$. This proves the required separation of the complete lines.

	Label the lower and upper left endpoints $a_j,b_j$ and the lower and upper right endpoints $c_j,d_j$. Then
	\begin{equation*}
		d(a_j,c_j)=d(b_j,d_j)=s,\qquad
		d(a_j,d_j)=d(b_j,c_j)=\sqrt{s^2+\ell^2}.
	\end{equation*}
	Thus \eqref{eq:cone-four-point-gap} holds with $q=4(\sqrt{s^2+\ell^2}-s)>0$.
\end{proof}

To control the measure of the corresponding families on $M$, it suffices to compare four distances. The next estimate avoids any assumption that measures on spaces of geodesics converge with the underlying metrics.

\begin{lemma}\label{lem:four-point-crofton}
	For $a,b,c,d\in M$, let $H([a,b],[c,d])$ be the family of oriented lines crossing both indicated geodesic segments. Then
	\begin{equation}\label{eq:four-point-crofton}
		m\bigl(H([a,b],[c,d])\bigr)
		\geq2\bigl|d(a,c)+d(b,d)-d(a,d)-d(b,c)\bigr|.
	\end{equation}
\end{lemma}

\begin{proof}
	A properly embedded line separates the plane into two components. Ignore the zero-measure pencils through the four endpoints, and for each remaining line denote the sides of an endpoint $x$ by $s_x\in\{0,1\}$. The indicator that this line crosses $[x,y]$ is $|s_x-s_y|$. Indeed, opposite sides force a crossing, whereas a transverse crossing between endpoints on the same side would force a second crossing. The latter is impossible for two distinct geodesic lines.

	The identity
	\begin{equation*}
		|s_a-s_c|+|s_b-s_d|-|s_a-s_d|-|s_b-s_c|
		=-2(s_a-s_b)(s_c-s_d)
	\end{equation*}
	does not depend on the labeling of the two sides. Integrate it and use \eqref{eq:crofton-segment}. The absolute value of the product on the right is the indicator of $H([a,b],[c,d])$, so
	\begin{equation*}
		4\bigl|d(a,c)+d(b,d)-d(a,d)-d(b,c)\bigr|
		\leq2m\bigl(H([a,b],[c,d])\bigr).
	\end{equation*}
	This is \eqref{eq:four-point-crofton}.
\end{proof}

\begin{lemma}\label{lem:transferred-families}
	Assume $0<E<\infty$. There is $q>0$, depending only on $E$, such that for every sequence $r_i\to\infty$ there are Borel families $H_i^1,H_i^2\subset\Gamma_{r_i}$ satisfying
	\begin{equation}\label{eq:positive-family-mass}
		\liminf_{i\to\infty}\frac{m(H_i^j)}{r_i}\ge q,
		\qquad j=1,2.
	\end{equation}
	For every fixed $T>2$, no pair in $H_i^1\times H_i^2$ has an intersection in $\overline B(p,(T+1)r_i)$ when $i$ is sufficiently large.
\end{lemma}

\begin{proof}
	By \autoref{thm:asymptotic-cone}, $(M,r_i^{-1}d,p)$ converges to $C_\Theta$ with $\Theta=2\pi+E>2\pi$. Take the four segments from \autoref{lem:cone-families}, choose approximating points for their eight endpoints, and join the corresponding endpoints in $M$ by geodesic segments. The limiting segments are unique, so the approximating segments converge to them. In particular, they lie in $B(p,r_i)$ for large $i$.

	For these indices, let $H_i^j$ consist of the oriented lines crossing both approximating segments of pair $j$; for the finitely many remaining indices, set $H_i^j=\varnothing$. These are Borel subsets of $\Gamma_{r_i}$. Equations \eqref{eq:four-point-crofton} and \eqref{eq:cone-four-point-gap}, applied to the converging endpoint distances, give \eqref{eq:positive-family-mass}.

	Fix $T>2$. If the last assertion failed, pass to a subsequence and choose a pair of lines meeting in $\overline B(p,(T+1)r_i)$. Anchor each line at its crossing with the first segment of its pair. In the rescaled metric, the second crossing occurs at a parameter of absolute value at most $2$, and the common intersection at a parameter of absolute value at most $T+2$. Compactness on bounded parameter intervals gives two complete limiting isometric lines. Each meets its pair of closed marked segments, and both contain the common limiting point. The marked points on each limiting line are distinct, since the two segments have positive separation. This contradicts \autoref{lem:cone-families}.
\end{proof}

The families in \autoref{lem:transferred-families} have product measure bounded below by a positive multiple of $r_i^2$. Cone convergence has excluded intersections in every fixed rescaled ball. It does not yet exclude intersections farther away. We now estimate precisely this remaining possibility.

\section{Distant intersections and rigidity}\label{sec:distant-intersections}

We continue to assume $E<\infty$. To complete the argument from \autoref{sec:separated-families}, we estimate the measure of pairs whose intersection is far from the sampling ball. For $T>2$, let
\begin{equation*}
N_r(T)=N_{\Gamma_r}\bigl(M\setminus B(p,(T+1)r)\bigr).
\end{equation*}
Thus $N_r(T)$ is the product measure of ordered pairs of lines in $\Gamma_r$ whose transverse intersection lies outside the indicated ball. Since Crofton measure scales as length, the appropriate normalization is $r^{-2}N_r(T)$.

\begin{proposition}\label{prop:intersection-tail}
	There is a constant $C_E<\infty$ such that, for every $T>2$,
	\begin{equation}\label{eq:intersection-tail}
		\limsup_{r\to\infty}\frac{N_r(T)}{r^2}\leq\frac{C_E}{T-1}.
	\end{equation}
\end{proposition}

\begin{proof}
	We first estimate the radial excess of a line, uniformly over all lines with an anchor in $\overline B(p,r)$. Set
	\begin{equation*}
		S_r(T)=\sup_{\gamma:\,d(p,\gamma(0))\leq r}
		\frac{d(p,\gamma(Tr))+d(p,\gamma(-Tr))-2Tr}{r}.
	\end{equation*}
	We claim that, for each fixed $T>2$,
	\begin{equation}\label{eq:uniform-radial-excess}
		\limsup_{r\to\infty}S_r(T)\leq\frac{1}{T}.
	\end{equation}
	Indeed, let $r_i\to\infty$ and let $\gamma_i$ be any sequence of such anchored lines. In the rescaled spaces $(M,r_i^{-1}d,p)$, the curves $t\mapsto\gamma_i(r_it)$ are isometric embeddings of $\mathbb R$, with anchors in the closed unit ball. By \autoref{thm:asymptotic-cone} and compactness on bounded parameter intervals, a subsequence converges to an isometric line $c:\mathbb R\to C_\Theta$, where $\Theta=2\pi+E$. Write $\rho(t)=d(0,c(t))$; then $\rho(0)\leq1$.

	If $c(a)=0$, isometry gives $\rho(t)=|t-a|$, with $|a|\leq1$. Otherwise, lift $c$ to the universal cover of the punctured cone and develop it into the Euclidean plane. Its image is an affine line avoiding the developed vertex: an affine line through that vertex would attain radius zero at a finite parameter. Consequently
	\begin{equation*}
		\rho(t)^2=(t-a)^2+b^2,\qquad b>0,\qquad a^2+b^2=\rho(0)^2\leq1.
	\end{equation*}
	In either case,
	\begin{equation*}
		\rho(T)^2+\rho(-T)^2=2T^2+2\rho(0)^2\leq2(T^2+1).
	\end{equation*}
	Hence Cauchy--Schwarz gives
	\begin{equation*}
		\rho(T)+\rho(-T)-2T
		\leq2\sqrt{T^2+1}-2T
		=\frac{2}{\sqrt{T^2+1}+T}\leq\frac1T.
	\end{equation*}
	A sequence violating \eqref{eq:uniform-radial-excess} would therefore produce a limiting line violating this bound. This proves the claim.

	For an anchored line put $h(t)=d(p,\gamma(t))$. Since $h(t)\leq r+|t|$, the part of the line outside $B(p,(T+1)r)$ has $|t|\geq Tr$. Apply \autoref{lem:angular-integral} to \eqref{eq:intersection-disintegration} on that exterior region and then use \autoref{lem:two-ended-excess}. The resulting intrinsic integral along each line satisfies
	\begin{equation*}
		\int_{\{t:\,h(t)\geq(T+1)r\}}(1-h'(t)^2)\,dt
		\leq2\bigl(h(Tr)+h(-Tr)-2Tr\bigr)
		\leq2r S_r(T).
	\end{equation*}
	The estimate holds for every allowable anchor, so it requires no measurable choice of anchors. Integrating gives
	\begin{equation*}
		N_r(T)\leq2C_0r m_r S_r(T).
	\end{equation*}
	The line-mass estimate \eqref{eq:line-mass-area} and \eqref{eq:uniform-area-bounds} give $m_r\leq C_E r$. Together with \eqref{eq:uniform-radial-excess}, this proves \eqref{eq:intersection-tail}.
\end{proof}

In \eqref{eq:intersection-tail}, $r$ tends to infinity with $T$ fixed. We may then let $T$ increase. This order of limits is exactly what is needed for the families constructed in \autoref{lem:transferred-families}.

\begin{proof}[Proof of \autoref{prop:probability-gap}]
	Let $\lambda=\limsup_{r\to\infty}P(r)$ and choose $r_i\to\infty$ with $P(r_i)\to\lambda$. Take the families $H_i^1,H_i^2$ from \autoref{lem:transferred-families}. By \eqref{eq:line-mass-area} and \eqref{eq:uniform-area-bounds}, there is a constant $B>0$, depending only on $E$, such that $m_r\le Br$ for all $r>0$.

	Fix $T>2$. For sufficiently large $i$, every transversely intersecting pair in $H_i^1\times H_i^2$ is counted by $N_{r_i}(T)$. All remaining pairs are among those in $\Gamma_{r_i}^2$ without a transverse intersection. Hence
	\begin{equation*}
		m(H_i^1)m(H_i^2)
		\le N_{r_i}(T)+(1-P(r_i))m_{r_i}^2.
	\end{equation*}
	Divide by $r_i^2$ and let $i\to\infty$. Equations \eqref{eq:positive-family-mass} and \eqref{eq:intersection-tail} give
	\begin{equation*}
		q^2\le\frac{C_E}{T-1}+(1-\lambda)B^2.
	\end{equation*}
	Letting $T\to\infty$ yields $\lambda\le1-q^2/B^2$. This proves the proposition with $\varepsilon=q^2/B^2>0$.
\end{proof}

\begin{proof}[Proof of \autoref{thm:intersection}]
	Suppose $\limsup_{r\to\infty}P(r)=1$. By \autoref{prop:probability-area}, the ratios $\operatorname{Area}(B(p,r_i))/r_i^2$ are bounded along some sequence $r_i\to\infty$, so \autoref{prop:integrability} gives $E<\infty$. \autoref{prop:probability-gap} excludes $E>0$. Thus $E=0$, and \autoref{prop:area-asymptotics} gives
	\begin{equation*}
		\lim_{r\to\infty}\frac{\operatorname{Area}(B(p,r))}{\pi r^2}=1.
	\end{equation*}
	The equality case of the area-growth theorem of Bangert and Emmerich \cite[Theorem~1]{BE2013} implies that $M$ is isometric to the Euclidean plane. Conversely, in the Euclidean plane the pairs of parallel lines have product Crofton measure zero, so $P(r)=1$ for every $r>0$.
\end{proof}

\section*{AI Assistance Disclosure}
This work was conducted with the assistance of AI. The central idea was generated by AI, and we subsequently organized and reformulated it.
AI was also used during the preparation of this manuscript
to assist with English-language writing and editing, as well as
the organization of relevant literature. 

\clearpage
\printbibliography

\end{document}